\documentclass[runningheads]{lncse}
\usepackage{amsmath}
\usepackage{amsfonts}
\usepackage{epsfig}
\usepackage{graphics}

\begin{document}

\title{Homogenization of the one-dimensional wave equation}
\author{ Thi Trang Nguyen \inst{1}, Michel Lenczner \inst{1} \and Matthieu
Brassart\inst{2} }

% your contribution title if the original one is too long use an abbreviated title (for running head):
\titlerunning{Homogenization of the one-dimensional wave equation}

% Use \authorrunning{Short Title} for an abbreviated version of the author list (for running head):
\authorrunning{Thi Trang Nguyen, Michel Lenczner, Matthieu Brassart.}

\institute{ FEMTO-ST, 26 Chemin de l'Epitaphe, 25000 Besan\c con,
France, {\tt thitrang.nguyen@femto-st.fr and
michel.lenczner@utbm.fr} \and Laboratoire de Math\'ematiques de
Besan\c con, 16 Route de Gray, 25030 Besan\c con, France, {\tt
matthieu.brassart@univ-fcomte.fr} } \maketitle

\begin{abstract}
We present a method for two-scale model derivation of the periodic
homogenization of the one-dimensional wave equation in a bounded
domain. It allows for analyzing the oscillations occurring on both
microscopic and macroscopic scales. The novelty reported here is on
the asymptotic behavior of high frequency waves and especially on
the boundary conditions of the homogenized equation. Numerical
simulations are reported.
\end{abstract}

\section{Introduction}

The paper is devoted to the periodic homogenization of the wave equation in
a one-dimensional open bounded domain where the time-independent
coefficients are $\varepsilon -$periodic with small period $\varepsilon >0$.
Corrector results for the low frequency waves have been published in \cite%
{brahim1992correctors,francfort1992oscillations}. These works were
not taking into account fast time oscillations, so the models
reflect only a part of the physical solution. In
\cite{brassart2009two}, an homogenized model has been developed to
cover the time and space oscillations occurring both at low and high
frequencies. Unfortunately, the boundary conditions of the
homogenized model was not found. Therefore, establishing the
boundary conditions of the homogenized model is critical and is the
goal of the present work which also extends
\cite{kader2000contributions}.

To this end, the wave equation is written under the form of a first order
formulation and the modulated two-scale transform $W_{k}^{\varepsilon }$\ is
applied to the solution $U^{\varepsilon }$\ as in \cite{brassart2009two}.
For $n\in
%TCIMACRO{\U{2115} }%
%BeginExpansion
\mathbb{N}
%EndExpansion
^{\ast }$\ and $k\in
%TCIMACRO{\U{211d} }%
%BeginExpansion
\mathbb{R}
%EndExpansion
,$\ the $n^{th}$\ eigenvalue $\lambda _{n}^{k}$\ of the Bloch wave problem
with $k$-quasi-periodic boundary conditions satisfies $\lambda
_{n}^{k}=\lambda _{n}^{-k}$, in addition $\lambda _{m}^{k}=\lambda _{n}^{k}$%
\ for $k\in
%TCIMACRO{\U{2124} }%
%BeginExpansion
\mathbb{Z}
%EndExpansion
/2$, so the corresponding waves are oscillating with the same frequency. The
homogenized model is thus derived for pairs of fibers $\{-k,k\}$ if $k\neq 0$
and for fiber $\left\{ 0\right\} $ otherwise\ which allows to derive the
expected boundary conditions. The weak limit of $\sum\nolimits_{\sigma \in
\left\{ -k,k\right\} }W_{\sigma }^{\varepsilon }U^{\varepsilon }$\ includes
low and high frequency waves, the former being solution of the homogenized
model derived in \cite{brahim1992correctors,francfort1992oscillations} and
the latter are associated to Bloch wave expansions. Numerical results
comparing solutions of the wave equation with solution of the two-scale
model for fixed $\varepsilon $ and $k$ are reported in the last section.

\section{The physical problem and elementary properties\label%
{Statement_problem}}

\textbf{The physical problem }We consider $I=\left( 0,T\right) \subset
\mathbb{R}^{+}$ a finite time interval and $\Omega =\left( 0,\alpha \right)
\subset \mathbb{R}^{+}$ a space interval, which boundary is denoted by $%
\partial \Omega $. Here, as usual $\varepsilon >0$ denotes a small parameter
intended to go to zero. Two functions $\left( a^{\varepsilon },\rho
^{\varepsilon }\right) $ are assumed to obey a prescribed profile $%
a^{\varepsilon }:=a\left( {\frac{x}{\varepsilon }}\right) $ and $\rho
^{\varepsilon }:=\rho \left( {\frac{x}{\varepsilon }}\right) $ where $\rho
\in L^{\infty }\left(
%TCIMACRO{\U{211d} }%
%BeginExpansion
\mathbb{R}
%EndExpansion
\right) $, $a\in W^{1,\infty }\left(
%TCIMACRO{\U{211d} }%
%BeginExpansion
\mathbb{R}
%EndExpansion
\right) $ are both $Y-$periodic where $Y=\left( 0,1\right) $. Moreover, they
are required to satisfy the standard uniform positivity and ellipticity
conditions, $0<\rho ^{0}\leq \rho \leq \rho ^{1}$ and $0<a^{0}\leq a\leq
a^{1},$ for some given strictly positive numbers $\rho ^{0}$, $\rho ^{1}$, $%
a^{0}$ and $a^{1}$. We consider $u^{\varepsilon }\left( t,x\right) $
solution to the wave equation with the source term $f^{\varepsilon }\in
L^{2}\left( I\times \Omega \right) $, initial conditions $%
(u_{0}^{\varepsilon },v_{0}^{\varepsilon })\in L^{2}\left( \Omega \right)
^{2}$ and homogeneous Dirichlet boundary conditions,%\vspace{-0.1cm}
\begin{equation}
\begin{array}{l}
\rho ^{\varepsilon }\partial _{tt}{u^{\varepsilon }}-{\partial _{x}}\left( {{%
a^{\varepsilon }}{\partial _{x}}u^{\varepsilon }}\right) ={f^{\varepsilon }}%
\,\ \text{in}\,\ I\times \Omega , \\
{u^{\varepsilon }}\left( {t=0,.}\right) =u_{0}^{\varepsilon }\,\ \text{and}%
\,\ {\partial _{t}}{u^{\varepsilon }}\left( {t=0,.}\right)
=v_{0}^{\varepsilon }\ \text{in}\,\ \Omega , \\
{u^{\varepsilon }}={0}\,\ \text{on}\,\ I\times \partial \Omega .%
\end{array}
\label{1D-wave-equation}
\end{equation}%
By setting: ${U^{\varepsilon }}:=({\sqrt{{a^{\varepsilon }}}{\partial _{x}}{%
u^{\varepsilon },}\sqrt{{\rho ^{\varepsilon }}}{\partial _{t}}{%
u^{\varepsilon })}},$ ${A^{\varepsilon }}=\left( {%
\begin{array}{cc}
0 & {\sqrt{{a^{\varepsilon }}}{\partial _{x}}\left( {\frac{1}{\sqrt{{\rho
^{\varepsilon }}}}.}\right) } \\
{\frac{1}{\sqrt{{\rho ^{\varepsilon }}}}{\partial _{x}}\left( {\sqrt{{%
a^{\varepsilon }}}.}\right) } & 0%
\end{array}%
}\right) ,\,$\ $U_{0}^{\varepsilon }:=({\sqrt{{a^{\varepsilon }}}{\partial
_{x}}u_{0}^{\varepsilon },\sqrt{{\rho ^{\varepsilon }}}v_{0}^{\varepsilon })}
$ and ${F^{\varepsilon }}:=(0,{{f^{\varepsilon }/}}\sqrt{{\rho ^{\varepsilon
}}})$, we reformulate the wave equation (\ref{1D-wave-equation}) as an
equivalent system,%\vspace{-0.1cm}%
\begin{equation*}
\left( {{\partial _{t}}-{A^{\varepsilon }}}\right) {U^{\varepsilon }}={%
F^{\varepsilon }}\text{ in }I\times \Omega ,{U^{\varepsilon }}\left( {t=0}%
\right) =U_{0}^{\varepsilon }\text{ in }\Omega \text{ and }{{{%
U_{2}^{\varepsilon }}}}={0}\text{ on }I\times \partial \Omega
\end{equation*}%
where ${{{U_{2}^{\varepsilon }}}}$ is the second component of $%
U^{\varepsilon }$. From now on, this system will be referred to as the
physical problem and taken in the distributional sense,%
\begin{equation}
\int\nolimits_{I\times \Omega }{{F^{\varepsilon }}\cdot {\Psi \,}+{%
U^{\varepsilon }}\cdot \left( {{\partial _{t}}-{A^{\varepsilon }}}\right)
\,\Psi dtdx}+\int\nolimits_{\Omega }{U_{0}^{\varepsilon }\cdot \Psi \left( {%
t=0}\right) \,dx}=0,  \label{1D-1st-weak-formulation}
\end{equation}%
for all the admissible test functions ${\Psi \in {H^{1}}{{\left( {I\times
\Omega }\right) }^{2}}}$ such that ${{\Psi \left( {t,.}\right) \in D\left( {{%
A^{\varepsilon }}}\right) }}$ for a.e. ${t\in I}$ where the domain $%
D(A^{\varepsilon }):=\{{\left( {\varphi ,\phi }\right) \in {L^{2}}\left(
\Omega \right) }^{2}|{{\sqrt{{a^{\varepsilon }}}\varphi \in {H^{1}}\left(
\Omega \right) ,}}$ ${{{\phi /\rho }\in {H_{0}^{1}}\left( \Omega \right) }\}}
$. As proved in \cite{brassart2009two}, the operator $iA^{\varepsilon }$
with the domain $D(A^{\varepsilon })$ is self-adjoint on $L^{2}(\Omega )^{2}$%
. We assume that the data are bounded ${\left\Vert {f^{\varepsilon }}%
\right\Vert _{{L^{2}}\left( {I\times \Omega }\right) }}+{\left\Vert {{%
\partial _{x}}u_{0}^{\varepsilon }}\right\Vert _{{L^{2}}\left( \Omega
\right) }}+{\left\Vert {v_{0}^{\varepsilon }}\right\Vert _{{L^{2}}\left(
\Omega \right) }}\leq c_{0}$, then $U^{\varepsilon }$ is uniformly bounded {%
in }${L^{2}}\left( {I\times \Omega }\right) .$

\textbf{Bloch waves }We introduce the dual $Y^{\ast }=\left( -\frac{1}{2},%
\frac{1}{2}\right) $ of $Y$. For any $k\in Y^{\ast }$, we define the space
of $k-$quasi-periodic functions $L_{k}^{2}:=\{u\in L_{loc}^{2}(\mathbb{R})$ $%
{|} $ $u(x+\ell )=u(x)e^{2i\pi k\ell }$ a.e. in $%
%TCIMACRO{\U{211d} }%
%BeginExpansion
\mathbb{R}
%EndExpansion
$ for all $\ell \in
%TCIMACRO{\U{2124} }%
%BeginExpansion
\mathbb{Z}
%EndExpansion
\}$ and set $H_{k}^{s}:=L_{k}^{2}\cap H_{loc}^{s}\left(
%TCIMACRO{\U{211d} }%
%BeginExpansion
\mathbb{R}
%EndExpansion
\right) $ for $s\geq 0.$ The periodic functions correspond to $k=0$. For a
given $k\in Y^{\ast }$, we denote by $(\lambda _{n}^{k},\phi _{n}^{k})_{n\in
%TCIMACRO{\U{2115} }%
%BeginExpansion
\mathbb{N}
%EndExpansion
^{\ast }}$ the Bloch wave eigenelements that are solution to%
\begin{equation*}
\mathcal{P}(k):-\partial _{y}\left( a\partial _{y}\phi _{n}^{k}\right)
=\lambda _{n}^{k}\rho \phi _{n}^{k}\text{ in }Y\text{ with }\phi _{n}^{k}\in
H_{k}^{2}(Y)\text{ and }\left\Vert \phi _{n}^{k}\right\Vert _{L^{2}\left(
Y\right) }=1.
\end{equation*}%
The asymptotic spectral problem $\mathcal{P}(k)$ is also restated as a first
order system by setting $A_{k}:=\left( {%
\begin{array}{cc}
0 & {\sqrt{a}\partial _{y}\left( \frac{1}{\sqrt{\rho }}.\right) } \\
\frac{1}{\sqrt{\rho }}{\partial _{y}\left( {\sqrt{a}.}\right) } & 0%
\end{array}%
}\right) $, $n_{A_{k}}=\frac{1}{\sqrt{\rho }}\left(
\begin{array}{cc}
0 & {\sqrt{a}n}_{Y} \\
{\sqrt{a}n}_{Y} & 0%
\end{array}%
\right) $ and $e_{n}^{k}:=\frac{1}{\sqrt{2}}\left(
\begin{array}{c}
-i{s_{n}/\sqrt{\lambda _{\left\vert n\right\vert }^{k}}}\sqrt{a}\partial
_{y}\left( {\phi _{\left\vert n\right\vert }^{k}}\right) \\
\sqrt{\rho }\phi _{\left\vert n\right\vert }^{k}%
\end{array}%
\right) $ where $s_{n}$ and $n_{Y}$ denote the sign of $n\in
%TCIMACRO{\U{2124} }%
%BeginExpansion
\mathbb{Z}
%EndExpansion
^{\ast }$ and the outer unit normal of $\partial Y$ respectively. As proved
in \cite{brassart2009two}, $iA_{k}$ is self-adjoint on the domain $D\left( {%
A_{k}}\right) :=\{{\left( {\varphi ,\phi }\right) \in L^{2}\left( Y\right)
^{2}|\sqrt{a}\varphi \in H_{k}^{1}\left( Y\right) ,}$ ${{\phi /}\sqrt{\rho }%
\in H_{k}^{1}\left( Y\right) }\subset L^{2}\left( Y\right) ^{2}\}.$ The
Bloch wave spectral problem $\mathcal{P(}k\mathcal{)}$ is equivalent to
finding pairs $\left( \mu _{n}^{k},e_{n}^{k}\right) $ indexed by $n\in
%TCIMACRO{\U{2124} }%
%BeginExpansion
\mathbb{Z}
%EndExpansion
^{\ast }$ solution to $\mathcal{Q(}k\mathcal{)}:A_{k}e_{n}^{k}=is_{n}\sqrt{%
\lambda _{\left\vert n\right\vert }^{k}}e_{n}^{k}$ \ in $Y$ with $%
e_{n}^{k}\in H_{k}^{1}\left( Y\right) ^{2}$. We pose $M_{n}^{k}:=\{m{\in {%
%TCIMACRO{\U{2124} }%
%BeginExpansion
\mathbb{Z}
%EndExpansion
}}^{\ast }{|}\lambda _{m}^{k}{=}\lambda _{n}^{k}$ and $s_{m}=s_{n}\}$ and
introduce the coefficients $b(k,n,m)=\int_{Y}\rho \phi _{\left\vert
n\right\vert }^{k}\cdot \phi _{\left\vert m\right\vert }^{k}dy$ and $%
c(k,n,m)=i{s_{n}/}\left( {2\sqrt{\lambda _{\left\vert n\right\vert }^{k}}}%
\right) \int_{Y}{\phi _{\left\vert n\right\vert }^{k}\cdot a{{\partial _{y}}}%
\phi _{\left\vert m\right\vert }^{k}}-a\partial _{y}{\phi _{\left\vert
n\right\vert }^{k}}\cdot {\phi _{\left\vert m\right\vert }^{k}}dy$ for $%
n,m\in M_{n}^{k}.$

\textbf{The modulated two-scale transform} Let us\textbf{\ }assume from now
that the domain $\Omega $ is the union of a finite number of entire cells of
size $\varepsilon $ or equivalently that the sequence $\varepsilon $ is
exactly $\varepsilon _{n}=\frac{\alpha }{n}$ for $n\in
%TCIMACRO{\U{2115} }%
%BeginExpansion
\mathbb{N}
%EndExpansion
^{\ast }$. For any $k\in Y^{\ast }$, we define $I^{k}=\left\{ -k,k\right\} $
if $k\neq 0$ and $I^{0}=\left\{ 0\right\} $. By choosing $\Lambda =\left(
0,1\right) $ as a time unit cell, we introduce the operator $%
W_{k}^{\varepsilon }:L^{2}\left( I\times \Omega \right) ^{2}\rightarrow
L^{2}\left( I\times \Lambda \times \Omega \times Y\right) ^{2}$ acting in
all time and space variables,\vspace{-0.1cm}
\begin{equation}
W_{k}^{\varepsilon }:=\left( 1-\sum\nolimits_{n\in
%TCIMACRO{\U{2124} }%
%BeginExpansion
\mathbb{Z}
%EndExpansion
^{\ast }}\Pi _{n}^{k}\right) S_{k}^{\varepsilon }+\sum\nolimits_{n\in
%TCIMACRO{\U{2124} }%
%BeginExpansion
\mathbb{Z}
%EndExpansion
^{\ast }}{{T^{\varepsilon \alpha _{n}^{k}}\Pi }}_{n}^{k}{S_{k}^{\varepsilon }%
}  \label{def-W}
\end{equation}%
where the time and space two-scale transforms ${{T^{\varepsilon \alpha
_{n}^{k}}}}$ and ${S_{k}^{\varepsilon }}$, and the orthogonal projector $\Pi
_{n}^{k}$ onto $e_{n}^{k}$ are defined in \cite{brassart2009two}, see pages
11,15 and 17, with $\alpha _{n}^{k}=2\pi /\sqrt{\lambda _{n}^{k}}$, and
where it is proved that,\vspace{-0.1cm}%
\begin{equation}
\left\Vert W_{k}^{\varepsilon }u\right\Vert _{L^{2}\left( I\times \Lambda
\times \Omega \times Y\right) }^{2}=\left\Vert u\right\Vert _{L^{2}\left(
I\times \Omega \right) }^{2}.  \label{two-scale-boundness}
\end{equation}%
We define $(\mathfrak{B}_{n}^{k}v)(t,x)=v(t,\frac{t}{\varepsilon \alpha
_{n}^{k}},x,\frac{x}{\varepsilon })$ the operator that operates on functions
$v(t,{\tau },x,y)$ defined in $I\times \mathbb{R\times }\Omega \times
\mathbb{R}$. The notation $O\left( \varepsilon \right) $ refers to numbers
or functions tending to zero when $\varepsilon \rightarrow 0$ in a sense
made precise in each case. The next Lemma shows that $\mathfrak{B}_{n}^{k}$
is an approximation of $T^{\varepsilon \alpha _{n}^{k}\ast
}S_{k}^{\varepsilon \ast }$ for a function which is periodic in $\tau $ and $%
k-$quasi-periodic in $y$, where $T^{\varepsilon \alpha _{n}^{k}\ast
}:L^{2}\left( I\times \Lambda \right) \rightarrow L^{2}\left( I\right) $ and
$S_{k}^{\varepsilon \ast }:L^{2}\left( \Omega \times Y\right) \rightarrow
L^{2}\left( \Omega \right) $ are adjoint of $T^{\varepsilon \alpha _{n}^{k}}$
and $S_{k}^{\varepsilon }$ respectively. \vspace{-0.15cm}%
\begin{lemma} \label{conver} Let $v\in C^{1}\left(
I\times \Lambda \times \Omega \times Y\right) $ a
periodic function in $\tau $ and $k-$quasi-periodic in $y$, then $%
T^{\varepsilon \alpha _{n}^{k}\ast }S_{k}^{\varepsilon \ast }v=\mathfrak{B}_{n}^{k}%
v+O\left( \varepsilon \right) $ in the $L^{2}\left( I \times\Omega
\right) $ sense. Consequently, for any sequence $u^{\varepsilon }$
bounded in $L^{2}\left( I\times \Omega \right) $ such that
$T^{\varepsilon \alpha
_{n}^{k}}S_{k}^{\varepsilon }u^{\varepsilon }$ converges to $u$ in $%
L^{2}(I\times \Lambda \times \Omega \times Y)$ weakly when $\varepsilon
\rightarrow 0$,%
\begin{equation}
\int_{I\times \Omega }u^{\varepsilon }\cdot \mathfrak{B}_{n}^{k}v\text{ }%
dtdx\rightarrow \int_{I\times \Lambda \times \Omega \times Y}u\cdot v\text{ }%
dtd\tau dxdy\text{ \ when }\varepsilon \rightarrow 0.  \label{convergence}
\end{equation}%
\end{lemma}\vspace{-0.15cm}Note that for $k=0$, the convergence (\ref%
{convergence}) regarding each variable corresponds to the definition of
two-scale convergence in \cite{allaire1992homogenization}. The proof is
carried out in three steps. First the explicit expression of $T^{\varepsilon
\alpha _{n}^{k}\ast }S_{k}^{\varepsilon \ast }v$ is derived, second the
approximation of $T^{\varepsilon \alpha _{n}^{k}\ast }S_{k}^{\varepsilon
\ast }v$ is deduced, finally the convergence (\ref{convergence}) follows.
For a function $v\left( t,\tau ,x,y\right) $ defined in $I\times \Lambda
\times \Omega \times Y,$ we observe that\vspace{-0.1cm}%
\begin{equation}
A^{\varepsilon }\mathfrak{B}_{n}^{k}v=\mathfrak{B}_{n}^{k}\left( \left(
\frac{A_{k}}{\varepsilon }+B\right) v\right) \text{ and }{{\partial _{t}}}%
\left( {\mathfrak{B}}_{n}^{k}v\right) ={\mathfrak{B}}_{n}^{k}\left( \left(
\frac{\partial _{\tau }}{\varepsilon \alpha _{n}^{k}}+\partial _{t}\right)
v\right) \text{,}  \label{derivative_t}
\end{equation}%
where the operator $B$ is defined as the result of the formal substitution
of $x-$derivatives by $y-$derivatives in $A_{k}$.

\section{Homogenized results and their proof\label{model}}

For $k\in Y^{\ast }$, we decompose\vspace{-0.1cm}
\begin{equation}
\frac{\alpha k}{\varepsilon }=h_{\varepsilon }^{k}+l_{\varepsilon }^{k}\text{
with }h_{\varepsilon }^{k}=\left[ \frac{\alpha k}{\varepsilon }\right] \text{
and }l_{\varepsilon }^{k}\in \left[ 0,1\right) ,  \label{epsilon_m}
\end{equation}%
and assume that the sequence $\varepsilon $ is varying in a set $%
E_{k}\subset \mathbb{R}^{+\ast }$ depending on $k$ so that\vspace{-0.1cm}%
\begin{equation}
l_{\varepsilon }^{k}\rightarrow l^{k}\text{ when }\varepsilon \rightarrow 0%
\text{ and }\varepsilon \in E_{k}\text{ with }l^{k}\in \left[ 0,1\right).
\label{l}
\end{equation}%
We note that for $k=0$, $h_{\varepsilon }^{k}=0,$ $l_{\varepsilon }^{k}=0$,
so $l^{k}=0$ and $E_{0}=\mathbb{R}^{+\ast }$. After extraction of a
subsequence, we introduce the weak limits of the relevant projections along $%
e_{n}^{k}$ for any $n\in {%
%TCIMACRO{\U{2124} }%
%BeginExpansion
\mathbb{Z}
%EndExpansion
}^{\ast }$,\vspace{-0.1cm}%
\begin{equation}
F_{n}^{k}:=\lim_{\varepsilon \rightarrow 0}\int\nolimits_{\Lambda \times Y}{T%
}^{\varepsilon \alpha _{n}^{k}}{S{_{k}^{\varepsilon }{F^{\varepsilon }}\cdot
{e^{2i\pi {s_{n}}\tau }}e_{n}^{k}dyd\tau }}\text{ and }U_{0,n}^{k}:=\lim_{%
\varepsilon \rightarrow 0}\int\nolimits_{Y}{S_{k}^{\varepsilon
}U_{0}^{\varepsilon }\cdot e_{n}^{k}dy}.  \label{data}
\end{equation}%
The next lemmas state the microscopic equation for each mode and the
corresponding macroscopic equation.\vspace{-0.1cm}%
\begin{lemma}
\label{lemma_micro}For $k\in Y^{\ast }$ and $n\in {%
%TCIMACRO{\U{2124} }%
%BeginExpansion
\mathbb{Z}
%EndExpansion
}^{\ast }$, let $U^{\varepsilon }$ be a bounded solution of (\ref%
{1D-1st-weak-formulation}), there exists at least a subsequence of $%
T^{\varepsilon \alpha _{n}^{k}}{{S}_{k}^{\varepsilon }U^{\varepsilon }}$
converging weakly towards a limit $U_{n}^{k}$ in $L^{2}(I\times \Lambda
\times \Omega \times Y)^{2}$ when $\varepsilon $ tends to zero. Then $%
U_{n}^{k}$ is a solution of the weak formulation of the microscopic equation%
\vspace{-0.2cm}\begin{equation} \left( \frac{{{\partial _{\tau
}}}}{\alpha _{n}^{k}}{-A}_{k}\right) U_{n}^{k}=0\text{ in }I\times
\Lambda \times \Omega \times Y \label{strong-form}
\end{equation}%
and is periodic in $\tau $ and $k-$quasi-periodic in $y$. Moreover,
it can be decomposed as
\begin{equation}
U_{n}^{k}\left( t,\tau ,x,y\right) ={\sum\limits_{p\in \mathbb{M}_{n}^{k}}}%
u_{p}^{k}\left( t,x\right) e^{{2{i\pi s}}_{p}\tau }e_{p}^{k}\left( y\right)
\text{ with }u_{p}^{k}\in L^{2}\left( I\times \Omega \right) .
\label{decompose_U}
\end{equation}%
\end{lemma}\vspace{-0.3cm}%
\begin{lemma}
\label{lemma_macro}For each $k\in Y^{\ast }$, $n\in {%
%TCIMACRO{\U{2124} }%
%BeginExpansion
\mathbb{Z}
%EndExpansion
}^{\ast }$, for each $\sigma \in I^{k}$ and $q\in M_{n}^{\sigma }$, the
macroscopic equation is stated by%
\begin{equation}
\begin{array}{l}
\sum\limits_{p\in {M}_{n}^{\sigma }}b\left( \sigma ,p,q\right) {\partial _{t}%
}u_{p}^{\sigma }-\sum\limits_{p\in {M}_{n}^{\sigma }}c\left( \sigma
,p,q\right) {\partial _{x}}u_{p}^{\sigma }=F_{q}^{\sigma }\,\ \text{in}\,\
I\times \Omega , \\
\sum\limits_{p\in {M}_{n}^{\sigma }}b\left( \sigma ,p,q\right) u_{p}^{\sigma
}\left( {t=0}\right) =U_{0,q}^{\sigma }\,\ \text{in}\,\ \Omega ,%
\end{array}
\label{macro}
\end{equation}%
with the boundary conditions in case where there exist $p\in
{M}_{n}^{k}$ such that $c\left( k,p,q\right) \neq 0$ and $\phi
_{\left\vert p\right\vert }^{k}(0)\ne0$
\vspace{-0.2cm}\begin{equation} \sum\limits_{\sigma \in
I^{k}}\sum\limits_{p\in {M}_{n}^{\sigma
}}u_{p}^{\sigma }\phi _{\left\vert p\right\vert }^{\sigma }\left( 0\right) {%
e^{sign\left( \sigma \right) 2i\pi \frac{l^{k}x}{\alpha }}}=0\,\ \text{on}%
\,\ I\times \partial \Omega .  \label{boundary}
\end{equation}%
\end{lemma}\vspace{-0.15cm}

The low frequency part $U_{H}^{0}$ relates to the weak limit in $L^{2}\left(
I\times \Omega \times Y\right) ^{2}$ of the kernel part of $%
S_{k}^{\varepsilon }$ in \ref{def-W}. It has been treated completely, in %
\cite{brahim1992correctors,brassart2009two}. Here, we focus on the
non-kernel part of $S_{k}^{\varepsilon }$, it relates to the high
frequency waves and microscopic and macroscopic scales. In order to
obtain the solution of the model, we analyze the asymptotic
behaviour of each mode through ${{T^{\varepsilon \alpha
_{n}^{k}}}S_{k}^{\varepsilon }}$ as in Lemma \ref{lemma_micro} and
Lemma \ref{lemma_macro}. Then the full solution is the sum of all
modes. We introduce the characteristic function $\chi _{0}\left(
k\right) =1$ if $k=0$ and $=0$ otherwise. The main Theorem states
as follows.\vspace{-0.1cm}%
\begin{theorem} \label{theorem}For a given \thinspace $k\in Y^{\ast
}$, let $U^{\varepsilon } $ be a solution of
(\ref{1D-1st-weak-formulation}) bounded in $L^{2}\left(
I\times \Omega \right) $, for $\varepsilon \in E_{k},$ as in (\ref{epsilon_m}%
, \ref{l}), the limit $G_{k}$ of any weakly converging extracted
subsequence of $\sum\limits_{\sigma \in I^{k}}W_{\sigma
}^{\varepsilon }U^{\varepsilon }$ in $L^{2}\left( I\times \Lambda
\times \Omega \times Y\right) ^{2}$ can be
decomposed as%
\vspace{-0.2cm}\begin{equation} G^{k}\left( t,\tau ,x,y\right) =\chi
_{0}\left( k\right) U_{H}^{0}\left( t,x,y\right)
+\sum\limits_{\sigma \in I^{k}}\sum\limits_{n\in
%TCIMACRO{\U{2124} }%
%BeginExpansion
\mathbb{Z}
%EndExpansion
^{\ast }}{u_{n}^{\sigma }\left( {t,x}\right) {e^{2i\pi {s_{n}}\tau }}%
e_{n}^{\sigma }\left( y\right) }  \label{decompose}
\end{equation}%
where $%
\left( u_{n}^{\sigma }\right) _{n,\sigma }$ are solutions of the
macroscopic equation (\ref{macro}, \ref{boundary}).
\end{theorem}\vspace{-0.15cm} Thus, it follows from (\ref{decompose}) that
the physical solution $U^{\varepsilon }$ is approximated by two-scale modes%
\vspace{-0.2cm}%
\begin{equation}
U^{\varepsilon }\left( t,x\right) \simeq \chi _{0}\left( k\right)
U_{H}^{k}\left( t,x,\frac{x}{\varepsilon }\right) +\sum\nolimits_{\sigma \in
I^{k}}\sum\nolimits_{n\in \mathbb{Z}^{\ast }}u_{n}^{\sigma }\left(
t,x\right) e^{is_{n}\sqrt{\lambda _{n}^{\sigma }}t/\varepsilon
}e_{n}^{\sigma }\left( \frac{x}{\varepsilon }\right) .
\label{physical_approximation}
\end{equation}%
The remain of this section provides the proofs of results.

\textbf{Proof of Lemma \ref{lemma_micro}}. %\begin{proof}
The test functions of the weak formulation (\ref{1D-1st-weak-formulation})
are chosen as $\Psi ^{\varepsilon }=\mathfrak{B}_{n}^{k}\Psi \left( {t,x}%
\right) $ for $k\in Y^{\ast }$, $n\in
%TCIMACRO{\U{2124} }%
%BeginExpansion
\mathbb{Z}
%EndExpansion
^{\ast }$ where $\Psi \in {C^{\infty } }\left( I\times
\Lambda\right.$$\left. \times \Omega
\times Y\right)^{2}$ is periodic in $\tau $ and $k-$quasi-periodic in $y$. From (%
\ref{derivative_t}) multiplied by $\varepsilon $, since ${\left( {\frac{{{%
\partial _{\tau }}}}{{\alpha _{n}^{k}}}-{A_{k}}}\right) \Psi }$ is periodic
in $\tau $ and $k-$quasi-periodic in $y$ and ${T^{\varepsilon \alpha
_{n}^{k}}S_{k}^{\varepsilon }{U^{\varepsilon }\rightarrow U}}_{n}^{k}$ in $%
L^{2}\left( I\times \Lambda \times \Omega \times Y\right) ^{2}$ weakly,
Lemma \ref{conver} allows to pass to the limit in the weak formulation, $%
\int\nolimits_{I\times \Lambda \times \Omega \times Y}{U}_{n}^{k}{\cdot
\left( \frac{{{\partial _{\tau }}}}{{\alpha _{n}^{k}}}{-}A_{k}\right) \Psi
dtd\tau dxdy}=0$. Using the assumption $U_{n}^{k}\in D\left( A_{k}\right)
\cap L^{2}\left( I\times \Omega \times Y;H^{1}\left( \Lambda \right) \right)
$ and applying an integration by parts,\vspace{-0.15cm}%
\begin{gather*}
\int\nolimits_{I\times \Lambda \times \Omega \times Y}\left( -\frac{{{%
\partial _{\tau }}}}{{\alpha _{n}^{k}}}+A_{k}\right) {U}_{n}^{k}{\cdot {\Psi
}dtd\tau dxdy+}\int_{{I\times \partial \Lambda \times \Omega \times Y}%
}U_{n}^{k}\cdot \Psi {dtd\tau dxdy} \\
-\int_{{I\times \Lambda \times \Omega \times \partial Y}}U_{n}^{k}\cdot
n_{A_{k}}\Psi {dtd\tau dxdy}=0.
\end{gather*}%
Then, choosing ${\Psi \in L}^{2}\left( I\times \Omega
;H_{0}^{1}\left( \Lambda \times Y\right) \right) $ comes the strong form (%
\ref{strong-form}). Since the product of a periodic function by a $k-$%
quasi-periodic function is $k-$quasi-periodic then $n_{A_{k}}\Psi $ is $k-$%
quasi-periodic in $y$. Therefore, $U_{n}^{k}$ is periodic in $\tau $ and $k-$%
quasi-periodic in $y.$ Moreover, (\ref{decompose_U}) is obtained, by
projection. %\end{proof}

\textbf{Proof of Lemma \ref{lemma_macro}} %\begin{proof}
For $k\in Y^{\ast }$, let $\left( \lambda _{p}^{\sigma },e_{p}^{\sigma
}\right) _{p\in M_{n}^{\sigma },\sigma \in I^{k}}$ be the Bloch eigenmodes
of the spectral equation $\mathcal{Q}\left( \sigma \right) $ corresponding
to the eigenvalue $\lambda _{n}^{k}$. We pose $\Psi ^{\varepsilon }\left(
t,x\right) =\sum\nolimits_{\sigma \in I^{k}}{\mathfrak{B}}_{n}^{k}\Psi
_{\varepsilon }^{\sigma }\in H^{1}\left( I\times \Omega \right) ^{2}$ as a
test function in the weak formulation (\ref{1D-1st-weak-formulation}) with
each $\Psi _{\varepsilon }^{\sigma }\left( t,\tau ,x,y\right)
=\sum\nolimits_{q\in M_{n}^{k}}\varphi _{q,\varepsilon }^{\sigma }\left(
t,x\right) e^{2i\pi s_{q}\tau }e_{q}^{\sigma }\left( y\right) $ where $%
\varphi _{q,\varepsilon }^{\sigma }\in H^{1}\left( I\times \Omega \right) $
and satisfies the boundary conditions

$\sum\nolimits_{\sigma \in I^{k},q\in M_{n}^{\sigma }}e^{2i\pi
s_{q}t/(\varepsilon \alpha _{q}^{\sigma })}\varphi _{q,\varepsilon }^{\sigma
}\left( t,x\right) \phi _{\left\vert q\right\vert }^{\sigma }\left( \frac{x}{%
\varepsilon }\right) =O\left( \varepsilon \right) $ on $I\times \partial
\Omega .$ Note that this condition is related to the second component of $%
\Psi ^{\varepsilon }$ only. Since $\alpha _{q}^{\sigma }=\alpha _{n}^{k}$
and $s_{q}=s_{n}$ for all $q\in M_{n}^{\sigma }$ and $\sigma \in I^{k}$, so $%
e^{2i\pi s_{q}t/(\varepsilon \alpha _{q}^{\sigma })}\neq 0$ can be
eliminated. Extracting a subsequence $\varepsilon \in E_{k}$, using the $%
\sigma -$quasi-periodicity of $\phi _{\left\vert q\right\vert }^{\sigma }$
and (\ref{epsilon_m},\ref{l}), $\varphi _{q,\varepsilon }^{\sigma }$
converges strongly to some $\varphi _{q}^{\sigma }$ in $H^{1}\left( I\times
\Omega \right) $, then the boundary conditions are\vspace{-0.15cm}%
\begin{equation}
\sum\nolimits_{\sigma \in I^{k}}\sum\nolimits_{q\in M_{n}^{\sigma }}\varphi
_{q}^{\sigma }\left( t,x\right) \phi _{\left\vert q\right\vert }^{\sigma
}\left( 0\right) e^{sign\left( \sigma \right) 2i\pi \frac{l^{k}x}{\alpha }%
}=0\,\text{on }I\times \partial \Omega .  \label{boundary-test}
\end{equation}%
Applying (\ref{derivative_t}) and since ${\left( {\frac{{{\partial _{\tau }}}%
}{{\alpha _{n}^{\sigma }}}-{{{A_{\sigma }}}}}\right) \Psi ^{\sigma }}=0$ for
$\sigma \in $ $I^{k}$, then in the weak formulation it remains\vspace{-0.15cm%
}%
\begin{equation*}
\sum\limits_{\sigma \in I^{k}}\int\nolimits_{I\times \Omega }{{%
F^{\varepsilon }\cdot \mathfrak{B}}}_{n}^{k}{{\Psi _{\varepsilon }^{\sigma
}+U^{\varepsilon }}\cdot \mathfrak{B}}_{n}^{k}{\left( {{\partial _{t}}-B}%
\right) \Psi _{\varepsilon }^{\sigma }dtdx}-\int\nolimits_{\Omega }{{%
U_{0}^{\varepsilon }}\cdot \mathfrak{B}}_{n}^{k}{\Psi _{\varepsilon
}^{\sigma }}\left( t=0\right) {dx=0}.
\end{equation*}%
Since ${\left( {{\partial _{t}}-B}\right) \Psi _{\varepsilon }^{\sigma }}$
is $\sigma -$quasi-periodic, so passing to the limit thanks to Lemma \ref%
{conver}, after using (\ref{data}) and replacing the decomposition of $%
U_{n}^{\sigma }$,\vspace{-0.15cm}%
\begin{gather*}
\sum\limits_{\sigma \in I^{k},\{p,q\}\in M_{n}^{\sigma }}\left(
\int\nolimits_{I\times \Omega }b\left( \sigma ,p,q\right) u_{p}^{\sigma
}\cdot {\partial }_{t}\varphi _{q}^{\sigma }-{c}\left( \sigma ,p,q\right) {{%
u_{p}^{\sigma }\cdot {{\partial _{x}}}}\varphi _{q}^{\sigma }-F_{q}^{\sigma }%
}\cdot {\varphi _{q}^{\sigma }\,dtdx}\right.  \\
\left. {-}\int_{\Omega }{U}_{0,q}^{\sigma }{\cdot \varphi _{q}^{\sigma }}%
\left( t=0\right) {dx}\right) {=0}\text{ for all }{\varphi _{q}^{\sigma }\in
}H^{1}\left( I\times \Omega \right) \text{ fulfilling (\ref{boundary-test})}.
\end{gather*}%
Moreover, if $u_{q}^{\sigma }{\in }H^{1}\left( I\times \Omega \right) $ then
it satisfies the strong form of the internal equations (\ref{macro}) for
each $\sigma \in I^{k}$, $q\in M_{n}^{\sigma }$ and the boundary conditions%
\vspace{-0.15cm}%
\begin{equation}
\sum\nolimits_{\sigma ,p,q}{{c\left( \sigma ,p,q\right) u_{p}^{\sigma }}}%
\overline{\varphi _{q}^{\sigma }}=0\text{ on }I\times \partial \Omega \text{
for }\varphi _{q}^{\sigma }\text{ satisfies (\ref{boundary-test})}.
\label{boudanry}
\end{equation}%
In order to find the boundary conditions of $\left( {{u_{p}^{\sigma }}}%
\right) _{\sigma ,p}$, we distinguish between the two cases $k\neq 0$ and $%
k=0$. First, for $k\neq 0$, $\lambda _{n}^{k}$ is simple so $%
M_{n}^{k}=\left\{ n\right\} $. Introducing $C=diag\left( c\left( \sigma
,n,n\right) \right) _{\sigma }$, $B=diag\left( b\left( \sigma ,n,n\right)
\right) _{\sigma }$, $U=\left( u_{n}^{\sigma }\right) _{\sigma }$, $F=\left(
F_{n}^{\sigma }\right) _{\sigma }$, $U_{0}=\left( U_{0,n}^{\sigma }\right)
_{\sigma }$, $\Psi =\left( \varphi _{n}^{\sigma }\right) _{\sigma }$, $\Phi
=\left( \phi _{\left\vert n\right\vert }^{\sigma }\left( 0\right)
e^{sign\left( \sigma \right) 2i\pi l^{k}x/\alpha }\right) _{\sigma }$,
Equation (\ref{macro}) states under matrix form\vspace{-0.15cm}%
\begin{equation}
B{\partial }_{t}U+C{\partial _{x}}U=F\text{ in }I\times \Omega \text{ and }%
BU\left( t=0\right) =U_{0}\text{ in }\Omega ,  \label{matrixform}
\end{equation}%
which boundary condition (\ref{boudanry}) is rewritten as $CU\left(
t,x\right) .\overline{\Psi }\left( t,x\right) =0$ on $I\times \partial
\Omega $ for all $\Psi $ such that $\overline{\Phi }(x).\overline{\Psi }%
(t,x)=0$ on $I\times \partial \Omega .$ Equivalently, $CU\left( t,x\right) $
is collinear with $\overline{\Phi }(x)$ yielding the boundary condition $%
u_{n}^{k}\phi _{\left\vert n\right\vert }^{k}\left( 0\right) {e^{2i\pi \frac{%
l^{k}x}{\alpha }}}$ ${+}$ $u_{n}^{-k}\phi _{\left\vert n\right\vert
}^{-k}\left( 0\right) {e^{-2i\pi \frac{l^{k}x}{\alpha }}}=0$ on $\ I\times
\partial \Omega $ after remarking that $c\left( k,n,n\right) \neq 0$ and $%
c\left( k,n,n\right) =-c\left( -k,n,n\right) $.

Second, for $k=0$, $\lambda _{n}^{0}$ is double $\lambda _{n}^{0}=\lambda
_{m}^{0}$ so $M_{n}^{k}=\left\{ n,m\right\} $. With $C=\left( c\left(
0,p,q\right) \right) _{p,q}$, $B=\left( b\left( 0,p,q\right) \right) _{p,q}$%
, $U=\left( u_{p}^{0}\right) _{p}$, $F=\left( F_{q}^{0}\right) _{q}$, $%
U_{0}=\left( U_{0,q}^{0}\right) _{q}$, $\Psi =\left( \varphi _{q}^{0}\right)
_{q}$, $\Phi =\left( \phi _{\left\vert q\right\vert }^{0}\left( 0\right)
\right) _{q}$, the matrix form is still stated as (\ref{matrixform}). Here,
the eigenvectors are chosen as real functions then $c\left( 0,p,p\right) =0.$
Since $c\left( 0,n,m\right) \neq 0$, so the boundary condition is $%
u_{n}^{0}\phi _{\left\vert n\right\vert }^{0}\left( 0\right) {+}%
u_{m}^{0}\phi _{\left\vert m\right\vert }^{0}\left( 0\right) =0\,\text{on}\
I\times \partial \Omega \text{.}$%\end{proof}

\textbf{Proof of Theorem} %\begin{proof}
For a given $k\in Y^{\ast }$, let $U^{\varepsilon }$ be solution of (\ref%
{1D-1st-weak-formulation}) which is bounded in $L^{2}(I\times \Omega )$, the
property (\ref{two-scale-boundness}) yields the boundness of $\left\Vert
W_{\sigma }^{\varepsilon }U^{\varepsilon }\right\Vert _{L^{2}\left( I\times
\Lambda \times \Omega \times Y\right) }$ for $\sigma \in I^{k}$. So there
exists $G^{k}\in L^{2}\left( I\times \Lambda \times \Omega \times Y\right)
^{2}$ such that, up to the extraction of a subsequence, $\sum\nolimits_{%
\sigma \in I^{k}}W_{\sigma }^{\varepsilon }U^{\varepsilon }$ tends weakly to
$G^{k}=\chi _{0}\left( k\right) U_{H}^{0}+\sum\nolimits_{\sigma \in
I^{k}}\sum\nolimits_{n\in
%TCIMACRO{\U{2124} }%
%BeginExpansion
\mathbb{Z}
%EndExpansion
^{\ast }}U_{n}^{k}$ in $L^{2}\left( I\times \Lambda \times \Omega \times
Y\right) ^{2}$. The high frequency part is based on the decomposition (\ref%
{decompose_U}) and Lemma \ref{lemma_macro}. %\end{proof}
\vspace{-0.15cm}%
\begin{remark} This method allows to complete the
homogenized model of the wave equation in \cite{brassart2009two} for
the one-dimensional case. Let $K\in
%TCIMACRO{\U{2115} }%
%BeginExpansion
\mathbb{N}
%EndExpansion
^{\ast }$, we decompose $\frac{\alpha }{\varepsilon K}=\left[ \frac{\alpha }{%
\varepsilon K}\right] +l_{\varepsilon }^{1}$ with $l_{\varepsilon
}^{1}\in \lbrack 0,1)$ and assume that the sequence $\varepsilon $
is varying in a set $E_{K}\subset
%TCIMACRO{\U{211d} }%
%BeginExpansion
\mathbb{R}
%EndExpansion
^{+\ast }$ depending on $K$ so that $l_{\varepsilon }^{1}\rightarrow
l^{1}$ when $\varepsilon \rightarrow 0$ with $l^{1}\in \lbrack
0,1)$. For any $k\in L_{K}^{\ast }$, defined in
\cite{brassart2009two}, we denote $p_{k}=kK\in
%TCIMACRO{\U{2115} }%
%BeginExpansion
\mathbb{N}
%EndExpansion
$, so $\frac{\alpha p_{k}}{\varepsilon K}=p_{k}\left[ \frac{\alpha }{%
\varepsilon K}\right] +p_{k}l_{\varepsilon }^{1}$ and
$p_{k}l_{\varepsilon }^{1}$ $\rightarrow l^{k}:=p_{k}l^{1}$ when
$\varepsilon \rightarrow 0$ with the same sequence of $\varepsilon
\in E_{K}$.
\end{remark}\vspace{-0.15cm}

\section{Numerical examples}

\begin{figure}[t]
\includegraphics[scale=.2]{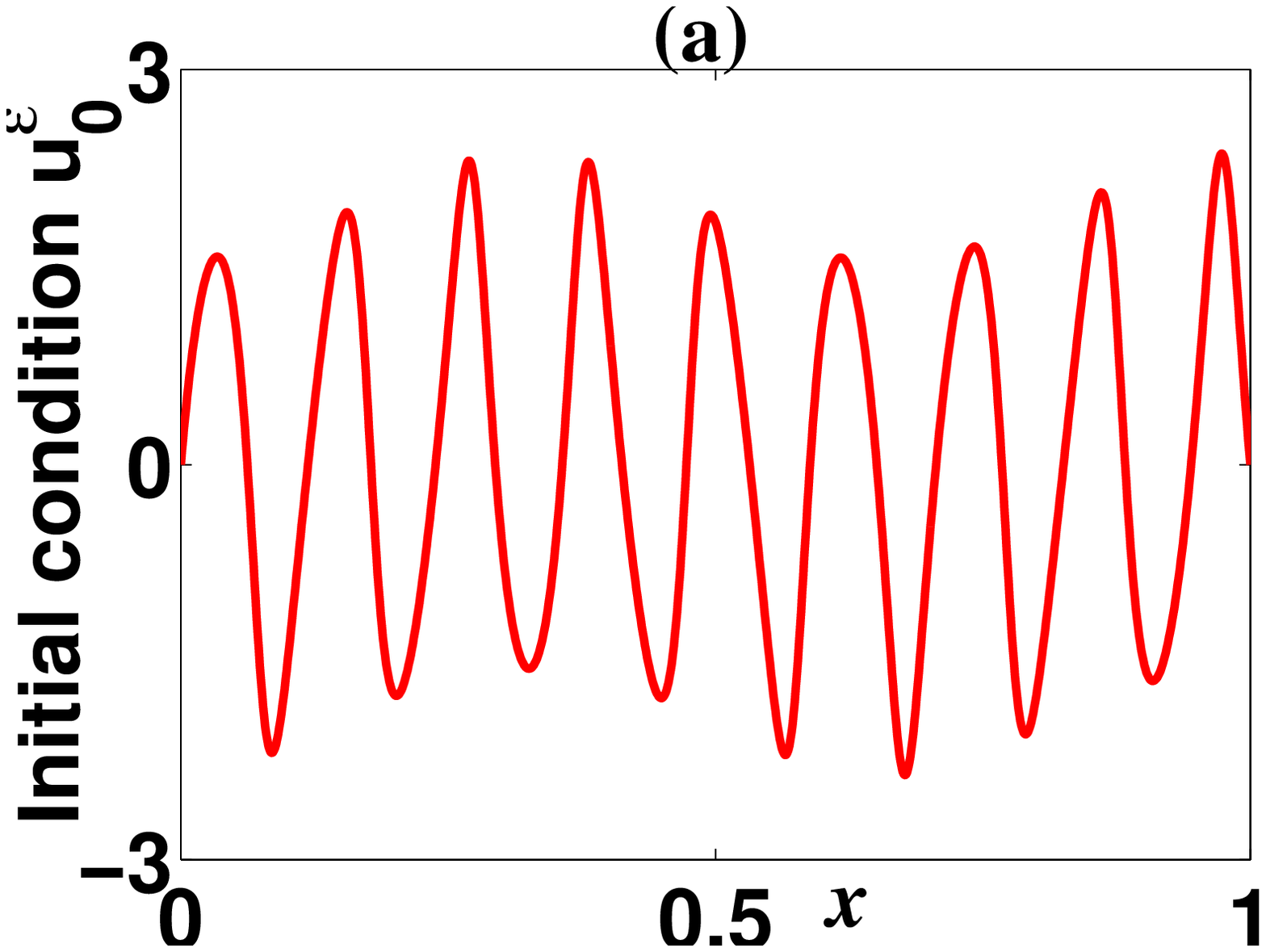}%
\includegraphics[scale=.2]{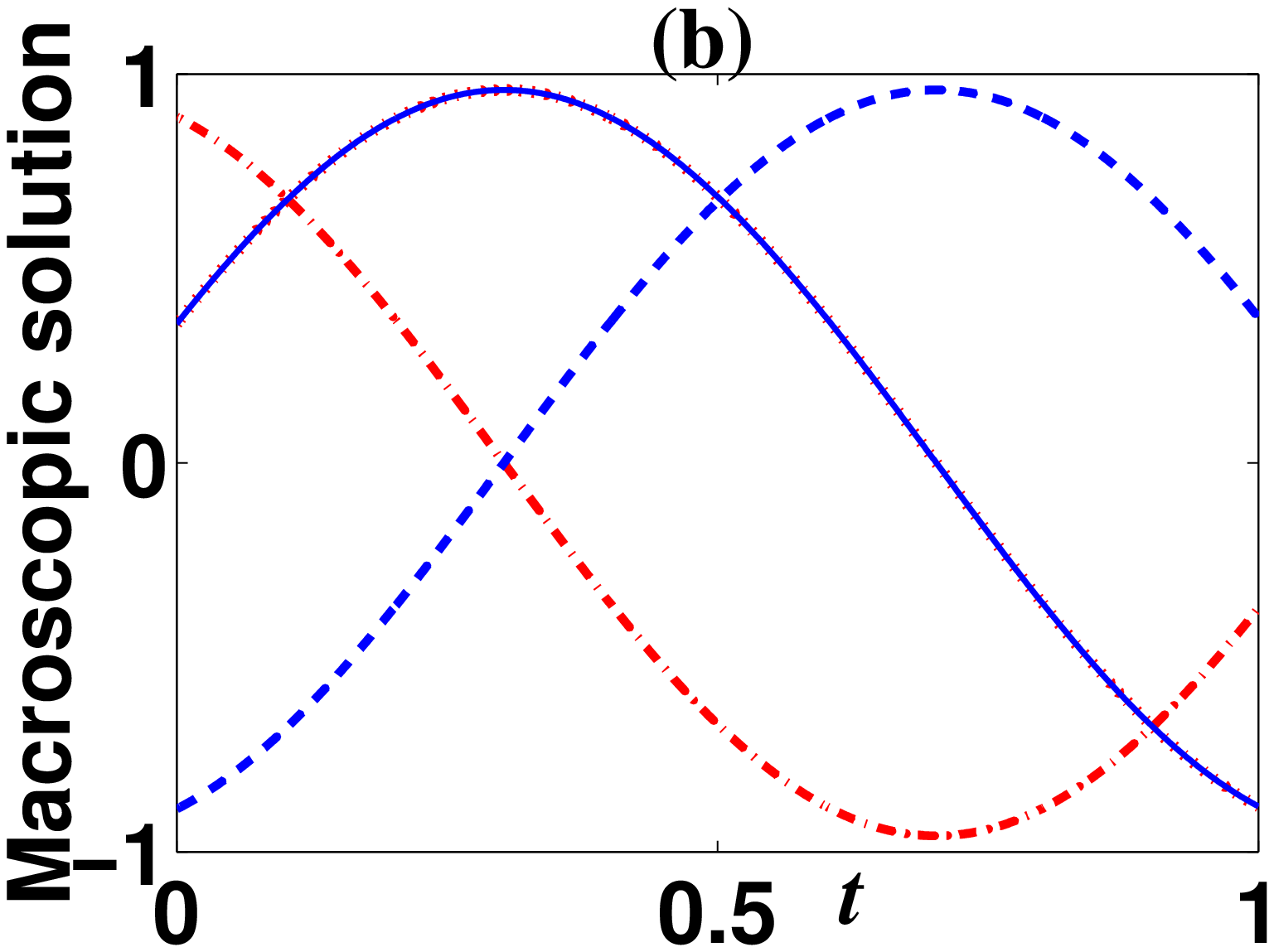}%
\includegraphics[scale=.2]{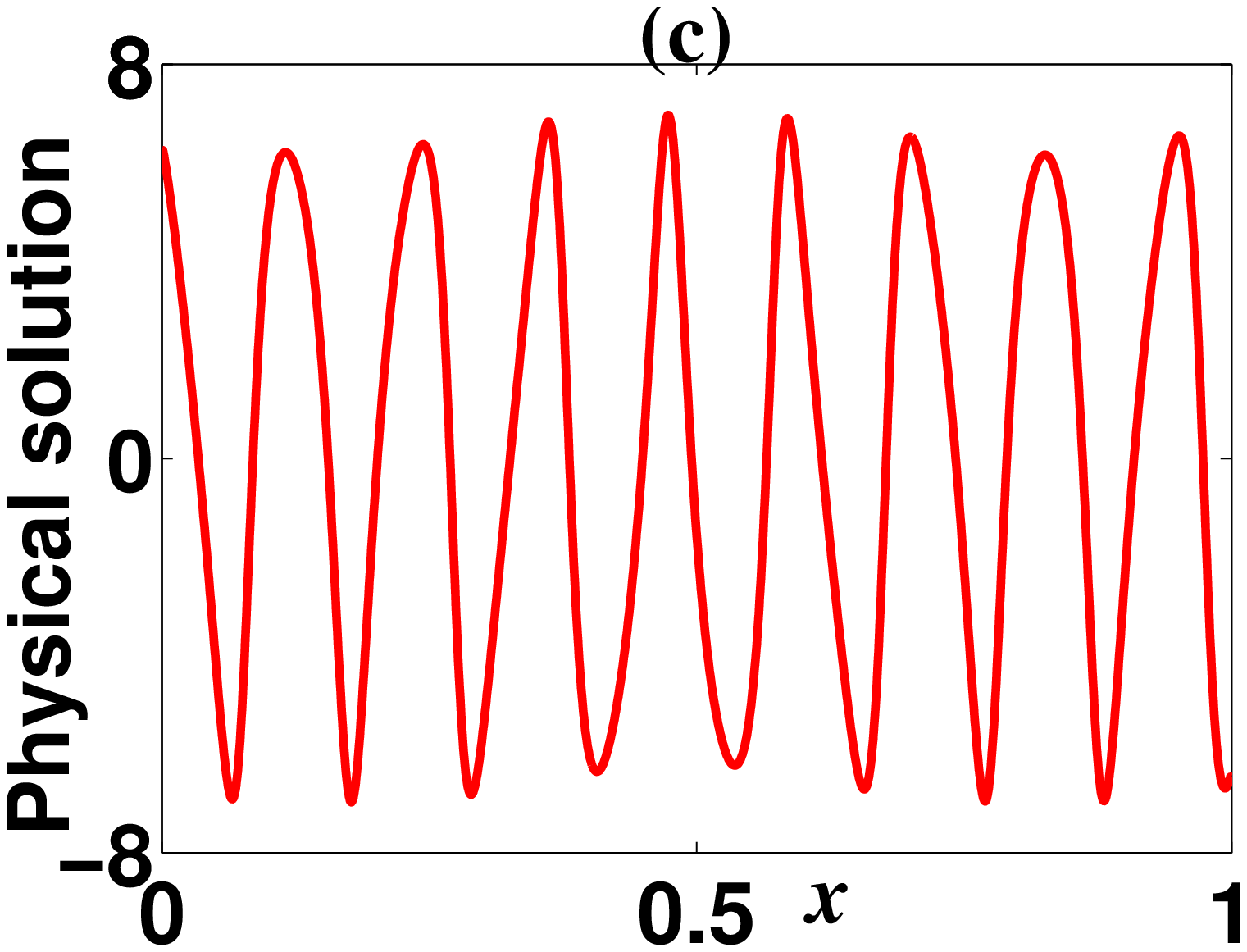}
\par
\includegraphics[scale=.2]{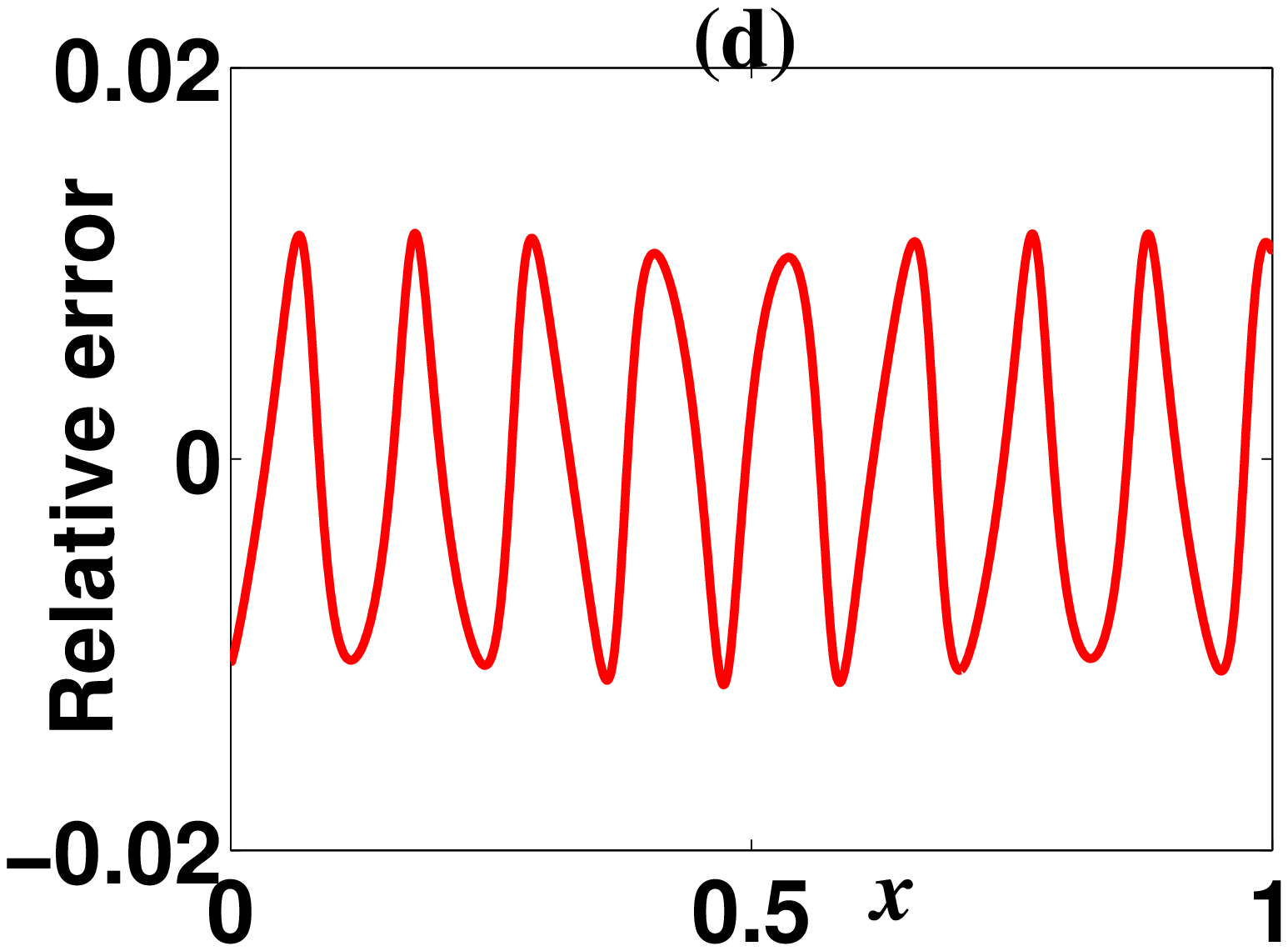}%
\includegraphics[scale=.2]{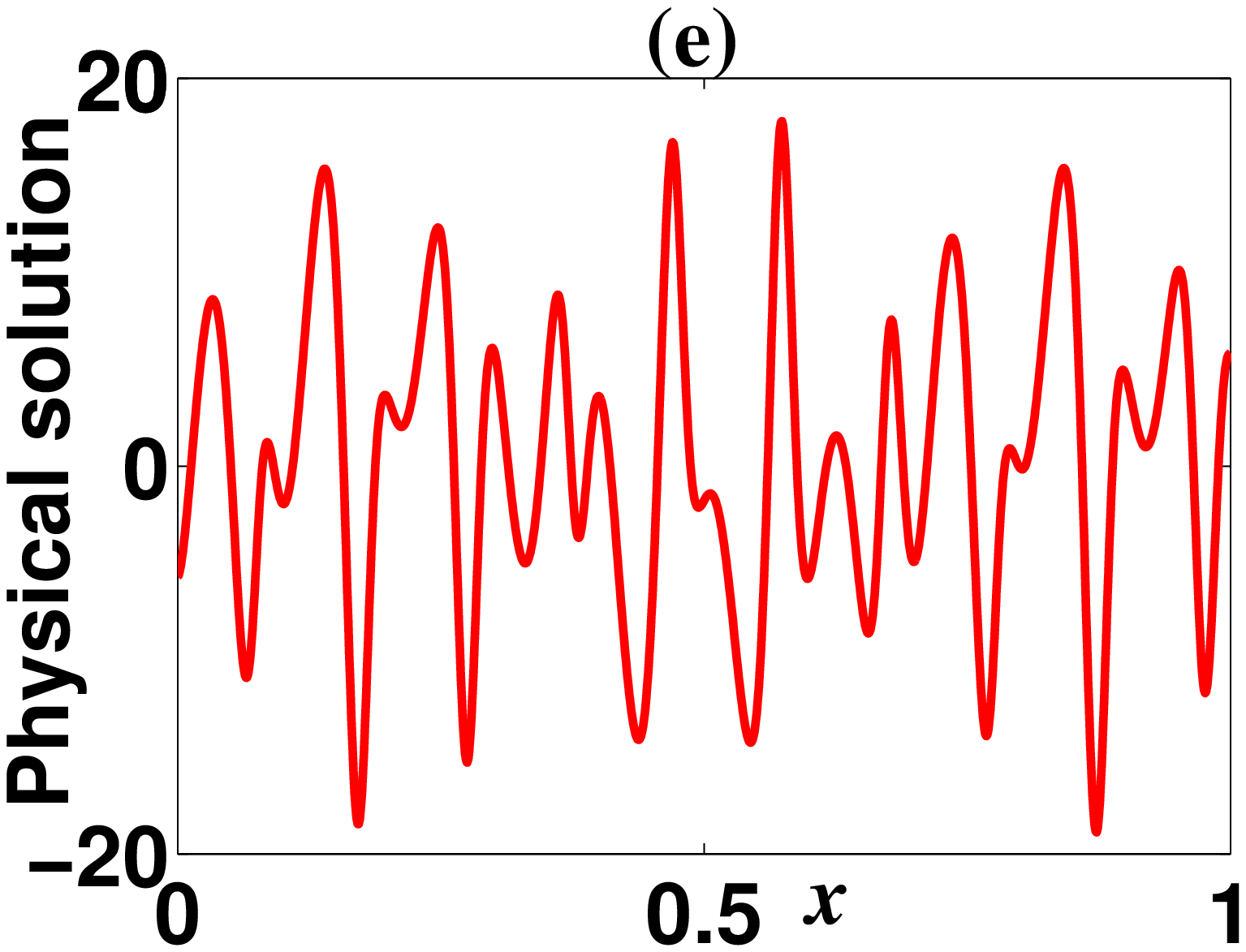}%
\includegraphics[scale=.2]{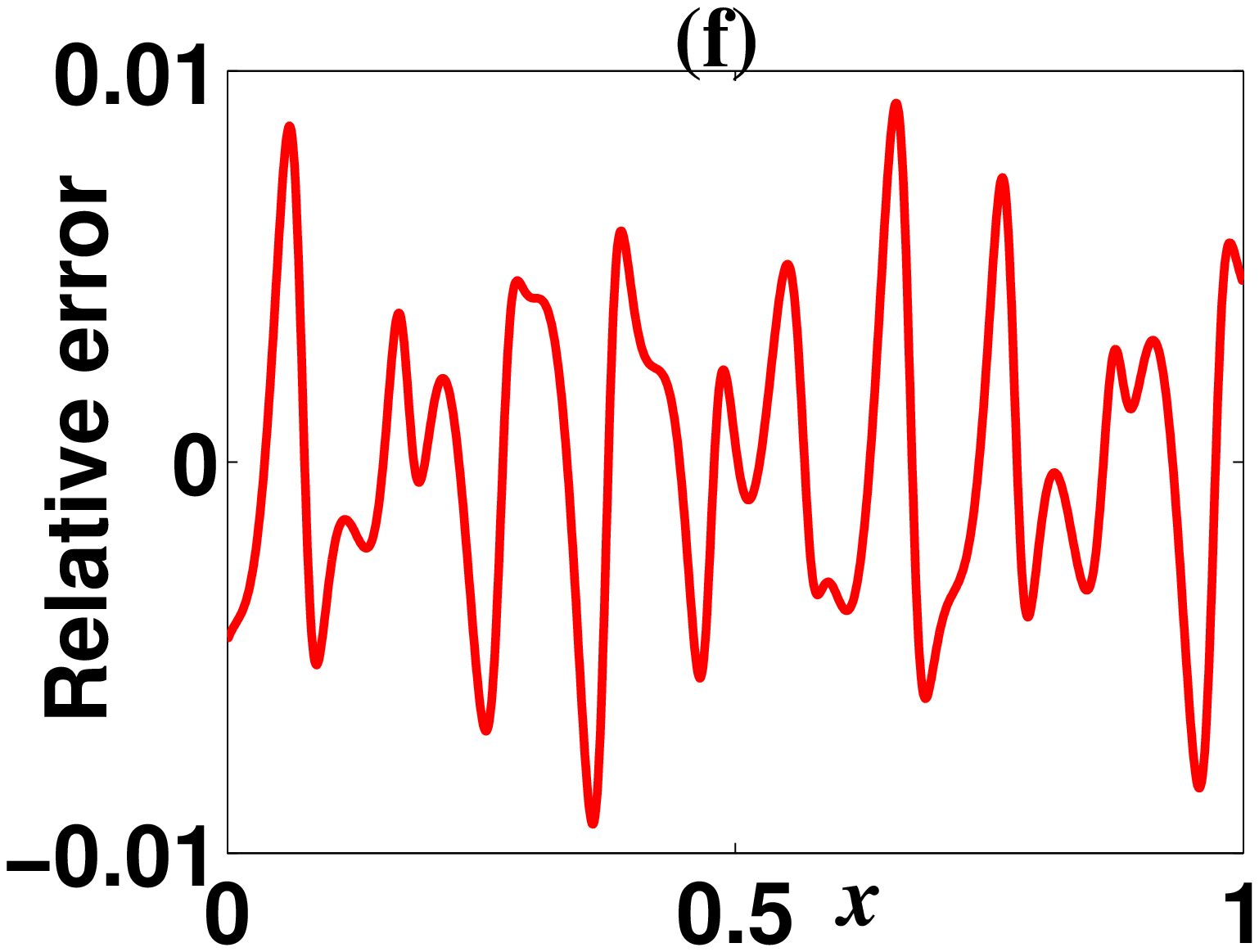} % Give a unique label
\caption{Numerical results}
\label{simulation}
\end{figure}

We report simulations regarding comparison of physical solution and its
approximation for $I=\left( 0,1\right) ,$ $\Omega =\left( 0,1\right) $, $%
\rho =1$, $a=\frac{1}{3}\left( \sin \left( 2\pi y\right) +2\right) $, $%
f^{\varepsilon }=0$,\thinspace\ $v_{0}^{\varepsilon }=0$, $\varepsilon =%
\frac{1}{10}$ and $k=0.16$. Since $k\neq 0$, so the approximation (\ref%
{physical_approximation}) comes\vspace{-0.2cm}%
\begin{equation}
U^{\varepsilon }\left( t,x\right) \simeq \sum\nolimits_{\sigma \in
I^{k}}\sum\nolimits_{n\in
%TCIMACRO{\U{2124} }%
%BeginExpansion
\mathbb{Z}
%EndExpansion
^{\ast }}u_{n}^{\sigma }\left( t,x\right) e^{is_{n}\sqrt{\lambda
_{n}^{\sigma }}t/\varepsilon }e_{n}^{\sigma }\left( \frac{x}{\varepsilon }%
\right) .  \label{num}
\end{equation}%
The validation of the approximation is based on the modal decomposition of
any solution $U^{\varepsilon }=\sum_{l\in
%TCIMACRO{\U{2124} }%
%BeginExpansion
\mathbb{Z}
%EndExpansion
^{\ast }}R_{l}^{\varepsilon }\left( t\right) V_{l}^{\varepsilon }\left(
x\right) $ where the modes $V_{l}^{\varepsilon }$ are built from the
solutions $v_{l}^{\varepsilon }$ of the spectral problem $\partial
_{x}\left( a^{\varepsilon }\partial _{x}v_{l}^{\varepsilon }\right) =\lambda
_{l}^{\varepsilon }v_{l}^{\varepsilon }$ in $\Omega $ with $%
v_{l}^{\varepsilon }=0$ on $\partial \Omega $.\ Moreover, in \cite%
{nguyen2013homogenization}, two-scale approximations of modes have been
derived on the form of linear combinations $\sum\nolimits_{\sigma \in
I^{k}}\theta _{n}^{\sigma }\left( x\right) \phi _{n}^{\sigma }\left( \frac{x%
}{\varepsilon }\right) $ of Bloch modes, so the initial conditions of the
physical problem are taken on the form \vspace{-0.2cm}%
\begin{equation}
u_{0}^{\varepsilon }\left( x\right) =\sum\nolimits_{n\in
%TCIMACRO{\U{2115} }%
%BeginExpansion
\mathbb{N}
%EndExpansion
^{\ast }}\sum\nolimits_{\sigma \in I^{k}}\theta _{n}^{\sigma }\left(
x\right) \phi _{n}^{\sigma }\left( \frac{x}{\varepsilon }\right) .
\label{initial condition}
\end{equation}%
Two simulations are reported, one for an initial condition $%
u_{0}^{\varepsilon }$ spanned by the pair of Bloch modes corresponding to $%
n=2$ when the other is spanned by three pairs $n\in \{2,3,4\}$. In the first
case, the first component of $U_{0}^{\varepsilon }$ approximates the first
component of a single eigenvector $V_{l}^{\varepsilon }$ approximated by (%
\ref{num}) where all coefficients $u_{n}^{\sigma }=0$ for $n\neq \pm 2$.
Fig. \ref{simulation} $(a)$ shows the initial condition $u_{0}^{\varepsilon
} $.\ Fig. \ref{simulation} $(b)$ presents the real part (solid line) and
the imaginary part (dashed-dotted line) of the macroscopic solution $%
u_{n}^{k}$ and also the real part (dotted line) and the imaginary part
(dashed line) of $u_{n}^{-k}$ at space step $x=0.699$ when Fig. \ref%
{simulation} $\left( c,d\right) $ plot the real part of the first component $%
U_{1}^{\varepsilon }$ of physical solution and the relative error vector of $%
U_{1}^{\varepsilon }$ with its approximation which $L^{2}(\Omega )$-norm is
equal to\textbf{\ }7e-3 at $t=0.466$\textbf{. }For the second case where $%
u_{n}^{\sigma }=0$ for $n\notin \{\pm 2,\pm 3,\pm 4\}$, the first component $%
U_{1}^{\varepsilon }$ and the relative error vector of $U_{1}^{\varepsilon }$
with its approximation which $L^{2}(\Omega )$-norm is 3.8e-3 are plotted in
Fig. \ref{simulation} $(e,f)$. Finally, for the two cases the $L^{2}(I)$%
-relative errors at $x=0.699$ on the first component are 8e-3 and 3.5e-3
respectively.
\bibliographystyle{vmams}
\bibliography{wave_reference}

\ifx\undefined\bysame
\newcommand{\bysame}{\leavevmode\hbox to3em{\hrulefill}\,}
\fi
\begin{thebibliography}{99}
\parskip1.0ex

\bibitem{allaire1992homogenization}
{\sc G.~Allaire}, {\em Homogenization and two-scale convergence}, SIAM Journal
  on Mathematical Analysis {\bf 23}:6 (1992), 1482--1518.

\bibitem{brahim1992correctors}
{\sc S.~Brahim-Otsmane, G.~Francfort, and F.~Murat}, {\em Correctors for the
  homogenization of the wave and heat equations}, Journal de math{\'e}matiques
  pures et appliqu{\'e}es {\bf 71}:3 (1992), 197--231.

\bibitem{brassart2009two}
{\sc M.~Brassart and M.~Lenczner}, {\em A two-scale model for the wave equation
  with oscillating coefficients and data}, Comptes Rendus Mathematique {\bf
  347}:23 (2009), 1439--1442.

\bibitem{francfort1992oscillations}
{\sc G.~A. Francfort and F.~Murat}, {\em Oscillations and energy densities in
  the wave equation}, Communications in partial differential equations {\bf
  17}:11-12 (1992), 1785--1865.

\bibitem{kader2000contributions}
{\sc M.~Kader}, {\em Contributions {\`a} la mod{\'e}lisation et contr{\^o}le
  des syst{\`e}mes intelligents distribu{\'e}s: Application au contr{\^o}le de
  vibrations d'une poutre}, Ph.D. thesis, Universit\'e de Franche-Comt\'e,
  France, 2000.

\bibitem{nguyen2013homogenization}
{\sc T.~T. Nguyen, M.~Lenczner, and M.~Brassart}, {\em Homogenization of the
  spectral equation in one-dimension}, arXiv preprint arXiv:1310.4064 (2013).

\end{thebibliography}

\end{document}